\documentclass{elsart3-r}
\usepackage{graphics}
\usepackage{amssymb}
\usepackage{amsfonts}
\usepackage{amsmath}
\usepackage[english,francais]{babel}
\usepackage{epsfig}
\usepackage[latin1]{inputenc}

\advance\textheight5mm

\newtheorem{e-proposition}[theorem]{Proposition}

\newtheorem{e-definition}[theorem]{Definition\rm}

\renewcommand{\theequation}{\arabic{equation}}

\def\og{\leavevmode\raise.3ex\hbox{$\scriptscriptstyle\langle\!\langle$~}}
\def\fg{\leavevmode\raise.3ex\hbox{~$\!\scriptscriptstyle\,\rangle\!\rangle$}}

\newcommand{\RR}{\field{R}}

\newcommand{\be}{\begin{equation}}
\newcommand{\ee}{\end{equation}}
\newcommand{\ba}{\begin{array}}
\newcommand{\ea}{\end{array}}
\newcommand{\bea}{\begin{eqnarray}}
\newcommand{\eea}{\end{eqnarray}}
\newcommand{\bee}{\begin{eqnarray*}}
\newcommand{\eee}{\end{eqnarray*}}

\catcode`@=11
\renewcommand\appendix{\bigskip {\noindent\Large \bf Appendix}\par
  \setcounter{section}{0}%
  \setcounter{subsection}{0}%
  \renewcommand\thesection{\@Alph\c@section}}
\catcode`@=12

\newenvironment{acknowledgement}{\noindent{\bf Acknowledgement.~}}{}

\def\RR{\mathbb{R}}
\def\R{{\mathbb R}}

\def\lim{\mathop{\rm lim}}

\def\exp{{\rm exp}}

\def\log{{\rm log}}

\def\un{{\bf 1}}

\def\fref#1{{\rm (\ref{#1})}}

\def\ds{\displaystyle}

\def\bs{\bigskip}

\newcommand{\dt}{\partial_t}

\newcommand{\grdx}{\nabla_x}
\newcommand{\eps}{\varepsilon}

\newcommand{\cint}[1]{\left\langle #1 \right\rangle}

\def\fref#1{{\rm (\ref{#1})}}

\begin{document}

\begin{frontmatter}
\selectlanguage{english}
\title{A boundary matching  micro/macro decomposition for kinetic equations}
\vspace{-2.6cm}
\selectlanguage{francais}
\title{Une nouvelle d\'ecomposition micro/macro adapt\'ee au bord pour les \'equations cin\'etiques}
\selectlanguage{english}
\author[authorlabel1]{Mohammed Lemou}
\ead{mohammed.lemou@univ-rennes1.fr}
\author[authorlabel2]{Florian M\'ehats}
\ead{florian.mehats@univ-rennes1.fr}
\address[authorlabel1]{IRMAR, CNRS \& Universit\'e de Rennes 1, Campus de beaulieu, 35042 Rennes, France}
\address[authorlabel2]{IRMAR \& Universit\'e de Rennes 1, Campus de beaulieu, 35042 Rennes, France}

\begin{abstract}
We introduce a new micro/macro decomposition of collisional  kinetic equations which naturally incorporates the exact space boundary conditions. The idea is to write the distribution fonction $f$ in all its domain as the sum of a Maxwellian adapted to the boundary (which is not the usual Maxwellian associated with $f$) and a reminder kinetic part. This Maxwellian is defined such that its 'incoming' velocity moments coincide with the 'incoming' velocity moments of the distribution function. Important consequences of this strategy are the following. i) No artificial boundary condition is needed in the micro/macro models and the exact boundary condition on $f$ is naturally transposed to the macro part of the model. ii) It provides a new class of the so-called 'Asymptotic preserving' (AP) numerical schemes: such schemes are consistent with the original kinetic equation for all fixed positive value of the Knudsen number $\eps$, and if $\eps \to 0 $ with fixed numerical parameters then these schemes degenerate into consistent numerical schemes for the various corresponding asymptotic fluid or diffusive models. Here, the strategy provides AP schemes not only inside the physical domain but also in the space boundary layers. We provide a numerical test in the case of a diffusion limit of the one-group transport equation, and show that our AP scheme recovers the boundary layer and a good approximation of the theoretical boundary value, which is usually  computed from to the so-called Chandrasekhar function.
\vskip 0.5\baselineskip
\selectlanguage{francais}
\noindent{\bf R\'esum\'e}
\vskip 0.5\baselineskip
\noindent
Nous introduisons une nouvelle d\'ecomposition micro/macro pour les \'equations cin\'etiques collisionnelles qui tient compte de fa\c{c}on exacte des conditions aux bords en espace. L'id\'ee est de d\'ecomposer la fonction de distribution $f$  en une Maxwellienne adapt\'ee au bord (qui n'est pas la Maxwellienne usuelle associ\'ee \`a $f$) et une partie cin\'etique restante. Cette nouvelle Maxwellienne est d\'efinie de sorte que ses moments  'entrants'  en vitesse soient les m\^emes que les moments 'entrants' de $f$.  Des cons\'equences importantes de cette strat\'egie sont les suivantes. i) Aucune condition artificielle n'est n\'ecessaire et les conditions de bord sur $f$ sont naturellement transpos\'ees sur la partie macroscopique de  la d\'ecomposition. ii) La strat\'egie produit une nouvelle classe de sch\'emas dits  'Asymptotic preserving' (AP) :  ces sch\'emas sont consistants avec le mod\`ele cin\'etique original pour toute valeur positive du nombre de Knudsen  $\eps$, et si $\eps \to 0 $ et les param\`etres num\'eriques sont fix\'es, alors ces sch\'emas d\'eg\'en\`erent en des sch\'emas qui sont consistants avec les divers mod\`eles asymptotiques associ\'es.  Ici, les sch\'emas obtenus reproduisent tr\`es bien les couches limites diffusives. Nous donnons un r\'esultat num\'erique dans le cas d'une limite de diffusion et montrons que la valeur au bord th\'eorique (calcul\'ee habituellement \`a partir de la fonction de Chandrasekhar) est bien approch\'ee num\'eriquement.
\end{abstract}
\end{frontmatter}

\selectlanguage{francais}
\section*{Version fran\c{c}aise abr\'eg\'ee}
Le but de ce travail est d'introduire une nouvelle d\'ecomposition micro/macro  pour les \'equations cin\'etiques, capable  par construction d'int\'egrer les conditions aux bords exactes du mod\`ele. Dans les approches usuelles, les mod\`ele macrosopiques (asymptotiques ou aux moments) utilisent g\'en\'eralement des moyennes en vitesse de la fonction $f$ ainsi que leurs flux, et les valeurs des ces moyennes au bord ne peuvent \^etre en g\'en\'eral d\'eduites de la condition au bord sur la fonction distribution $f$. En effet, si la condition au bord en espace est de type Dirichlet sur $f$, alors la valeur de $f$ au bord ne peut \^etre impos\'ee que pour les vitesses entrantes. Par cons\'equent les moyennes (sur tout le domaine des vitesses) de $f$ ne sont pas donn\'ees a priori, et des calculs suppl\'ementaires (probl\`eme de Milne) et/ou des conditions artificielles sont en g\'en\'eral n\'ecessaires pour bien poser le mod\`ele macroscopique. 

La strat\'egie pr\'esent\'ee ici permet de r\'esoudre ce probl\`eme et constitue une base pour le d\'eveloppement des sch\'emas dit  'Asymptotic Preserving' (AP) qui ont en plus  la sp\'ecifit\'e importante de pouvoir int\'egrer les conditions  aux bords de fa\c{c}on exacte. Avant de d\'etailler la strat\'egie dans les sections suivantes, nous rappelons bri\`evement la probl\'ematique des sch\'emas multi-\'echelles dans ce contexte. Quand le nombre de Knudsen $\eps$  devient petit,  une raideur appara\^{\i}t  dans ces mod\`eles  qui fait en sorte que leurs simulations  num\'eriques  par  des sch\'emas explicites en temps  deviennent rapidement inaccessibles (contraintes de type $\Delta t=O(\eps)$ ou $O(\eps^2)$). Le probl\`eme est alors de construire des sch\'emas num\'eriques pour l'\'equation cin\'etique consid\'er\'ee qui s'affranchissenent  de ces contraintes et d\'eg\'en\`erent en des discr\'etisations consistantes avec les  mod\`eles asymptotiques quand $\eps \to 0$.  Plusieurs travaux ont \'et\'e effectu\'es \`a ce propos, voir par exemple  \cite{JPT,Klar1,Klar2} pour les limites de diffusion. Une bibliographie plus compl\`ete sera donn\'ee dans une version d\'etaill\'ee de ce travail \cite{Lem-Meh-AP1}.
Une approche g\'en\'erale a \'et\'e introduite dans \cite{BLM,LM-AP1,Lem-note1}, permettant de construire des sch\'emas num\'eriques qui sont consistants avec les mod\`eles cin\'etiques \'etudi\'es pour toutes les valeurs de $\eps>0$, et d\'eg\'en\`erent quand $\eps \to 0$ en des discr\'etisations consistantes avec les syst\`emes d'Euler,  de Navier-Stokes compressibles, ou avec la limite de diffusion.  Elle est bas\'ee sur la d\'ecomposition micro-macro qui permet d'\'ecrire de mani\`ere \'equivalente l'\'equation cin\'etique de d\'epart sous la forme d'un syst\`eme couplant une \'equation fluide ou de diffusion avec une \'equation sur la partie cin\'etique restante. 

Le but ici est de compl\'eter cette s\'erie de travaux sur les sch\'emas micro/macro, en d\'eveloppant une nouvelle strat\'egie qui tient compte des conditions de bord en espace de fa\c{c}on exacte. Elle consiste \`a choisir une partie macrosopique 'int\'egr\'ee sur les vitesses entrantes', qui est compatible avec les conditions aux bords, et de la coupler avec la partie cin\'etique restante. En plus d'int\'egrer naturellement les conditions aux bords, cette strat\'egie constitue une base pour la construction des sch\'emas num\'eriques qui ont  toutes les bonnes propri\'et\'es des sch\'emas micro/macro d\'ecrits ci-dessus. Cette note donne les grandes lignes de la strat\'egie qui sera d\'evelopp\'ee en d\'etail dans \cite{Lem-Meh-AP1}. Nous donnons aussi les r\'esultats d'un premier test num\'erique dans le cas d'une limite diffusive, qui montre que le sch\'ema construit est bien AP et permet de retrouver la bonne couche limite au bord. Des exp\'eriences num\'eriques plus compl\`etes seront effectu\'ees dans \cite{Lem-Meh-AP1}.
\selectlanguage{english}
\section{Boundary matching micro/macro formulations of kinetic equations: the general problem}

It is a usual challenge to design efficient numerical schemes for kinetic equations which are consistent with the kinetic model for all positive value of the Knudsen number $\eps$, and degenerate into consistent schemes with the asymptotic models (compressible Euler and Navier-Stokes, diffusion, etc) when $\eps \to 0$. On this subject, many works can be quoted, see for instance \cite{JPT,Klar1,Klar2} and a more complete bibliography in the forthcoming paper \cite{Lem-Meh-AP1}. Recently  \cite{BLM,LM-AP1,Lem-note1}, a general strategy based on the micro/macro decomposition was introduced. It consists in writing the kinetic equation as a system coupling a macroscopic part (say a Maxwellian) and a reminder kinetic part.  The main advantage of this method is its robustness and easy adaptability.  
However, the general problem of  transfering boundary conditions from the distribution function $f$ to the two parts of the micro/macro systems has not been  solved and artificial boundary conditions were needed.

Our aim is to develop a new micro/macro decomposition of  collisional  kinetic equations which naturally incorporates the exact  space boundary conditions (BC). In fact, it is well known that during the derivation of a macroscopic model from kinetic equation (via moments or Chapman-Enskog approach), the exact boundary conditions on $f$ are generally lost. Indeed, if one takes moments of the kinetic equation, then the values of the fluxes at the boundary are required. But, in general (for inflow BC for instance), the distribution function $f$ is  known on the boundary for only incoming velocities and, therefore, their moments are not prescribed at the boundary. 

The idea is to decompose the distribution fonction $f$ in its domain as the sum of a  Maxwellian part adapted to the boundary (which is not the usual Maxwellian associated with $f$) and a reminder kinetic part. This Maxwellian is defined on the whole domain such that its 'incoming' velocity moments coincide with the 'incoming' velocity moments of the distribution function.  Important consequences of this strategy are the following. i) No artificial boundary condition is needed in the micro/macro models and the exact boundary conditions on $f$ are naturally shared by the macro part and the kinetic part. Note that the traditional approaches cannot provide exact boundary conditions for the asymptotic and kinetic  models separately, and extra calculations (Milne problems for instance) or artificial boundary conditions are necessary in general. ii) It provides a new  class of the so-called 'Asymptotic preserving' (AP) numerical schemes: such schemes are consistent with the original kinetic equation for all fixed positive value of the Knudsen number $\eps$, and when $\eps \to 0 $ with fixed numerical parameters, these schemes degenerate into consistent numerical schemes  for the various corresponding asymptotic models.  A further fundamental property of the strategy developed here is that it provides AP schemes  not only inside the physical domain but also in the space boundary layers with exact boundary conditions. We shall refer to this property as 'Boundary  Asymptotic Preserving' (BAP). At the end of this Note, we provide a numerical test which illustrates this property. We emphasize that our scheme does not use any artificial boundary condition: it provides at the diffusion limit a good approximation of the theoretical boundary value derived from the Chandrasekhar function, without injecting this value in our scheme.
This Note summarizes the main lines of the strategy and announces a more complete presentation \cite{Lem-Meh-AP1}, where numerical discretizations of this strategy will be constructed and implemented with a large variety of numerical tests.

\vspace*{-5mm}
\section{The case of the diffusion scaling}
Let $\Omega$ be a domain in the position (physical) space  $\RR^d$, with boundary $\partial \Omega$, and $V$ be a domain in the velocity space $\RR^d$, endowed with a measure $d\mu$.  We consider a function $\omega: \Omega\times V \rightarrow \RR$ such that for all $x\in \partial \Omega$ and $v\in V$, $\omega(x,v)$ has the same sign as $n(x)\cdot v$, where $n(x)$ is the outgoing normal vector to $\partial \Omega$ at $x$. Let us give explicit examples of $\omega(x,v)$ for specific geometries. For a ball centered at the origin, one can take $\omega(x,v)=x\cdot v$. For a half plane $x_{1}>0$, one can choose $\omega(x,v)=(-v_1, 0, ...,0)$.  In dimension one, for the interval $(0,1)$, one can take $\omega(x,v)=(2x-1)v$. We then consider the following transport equation in a diffusive scaling
\begin{equation}  \label{eq-transp}
\eps \dt f + v\cdot \nabla_x f = \frac{1}{\eps}L f ,   \qquad t > 0,  \qquad
  (x,v) \in \Omega \times V, \qquad f|_{t=0}=f_{init} 
\end{equation}
where $f$ is the distribution function of the particles that
depends on time $t>0$, on position $x\in\Omega $, and on velocity
$v\in V$.  The linear operator $L$ acts on the velocity dependence of $f$ and describes the interactions of particles with the medium.
We assume that there exists a positive equilibrium function $\mathcal E= \mathcal E(v)$ with $\cint{\mathcal E}:=\int_{V} \mathcal E d\mu =1$ and satisfying $\cint{v\mathcal E}=0$, and that  the collision operator $L$ is non-positive and self-adjoint in ${\rm L}^2(V,\mathcal E^{-1}\, d\mu)$;
 with null space and  image given by
 ${\mathcal N}(L)=\mbox{Span}\lbrace\mathcal E\rbrace=\lbrace f=\rho \mathcal E, \text{where } \rho:=\cint{f}:=\int_V fd\mu\rbrace, \quad 
{\mathcal R}(L) = ({\mathcal N}(L))^\bot = \lbrace f \text{ such that } \cint{f}=0 \rbrace.
$
When $\eps$ goes to 0 in~(\ref{eq-transp}), it is easily seen that $f$
converges to an equilibrium state $f_0=\rho_0(t,x) \mathcal E(v)$. The
diffusion limit is the equation satisfied by the density $\rho_0$ and is classically given by
\begin{equation}  \label{eq-diff}
\dt \rho_0 + \nabla_x \cdot (\kappa \nabla_x \rho_0) = 0, \qquad  \mbox{with}\qquad \kappa=\cint{vL^{-1}(v\mathcal E)}.
\end{equation}

Following \cite{LM-AP1} and \cite{Lem-note1},   equation \fref{eq-transp} can be written in a micro-macro equivalent form via the decomposition $f= \rho \mathcal E + g, \ \mbox{with} \  \rho(t,x)= \cint{f}$ and $g=f-\rho \mathcal E$ (g is not necessarily small as in \cite{LM-AP1}).  Let $\Pi$ be the orthogonal projector in $L^2(\mathcal E^{-1}d\mu)$ onto the nullspace of $L$:
$\Pi \phi= \cint{\phi}  \mathcal E$.
Then, inserting this decomposition  into the kinetic equation and applying  $\Pi$ and $I -\Pi$ successively, one gets after direct computations
\be
\label{kinetic-diffusion}
\ds \dt \rho + \frac{1}{\eps}  \grdx \cdot
  \cint{ v g } =  0, \ \ \ \ \ \  
   \dt g + \frac{1}{\eps} (I - \Pi)(v \cdot \grdx g) = \frac{1}{\eps^2}
    \big[ L g - \eps  \mathcal Ev\cdot \grdx \rho) \big] .
\ee

We now emphasize that in this formulation, the space boundary condition on $\rho$ and $g$ are not known because they cannot be inferred from the boundary conditions on $f$ in general.  Indeed, in the typical case of incoming boundary conditions, we impose 
\be
\label{BC-f}
f(t,x,v)= f_{b}(t,x,v), \qquad \forall t>0, \qquad \forall (x,v) \in \partial \Omega \times V\ \mbox{such that} \ \omega(x, v) < 0.
\ee
It is therefore clear that the values of $\rho(t,x) =\int_{V} f( t,x,v)d\mu$  cannot be determined a priori on the boundary $\partial \Omega$, since $f$ is only known for incoming velocities $v$, i.e. such that $\omega(x,v) < 0.$
Our aim in this work is to develop a new micro/macro decomposition  which is able to incorporate the boundary conditions in an exact way. To this purpose, let us introduce some few further notations:
\be
\label{VOpm}
\begin{array}{ll}
V_{-}(x)= \{v\in V, \ \omega(x,v)<0 \},  &\qquad \    \ V_{+}(x)= V\backslash V_{-}(x).
\end{array}
\ee
The idea is now the following. Instead of looking for an equation on $\rho$ as usual, we seek  an equation on the following 'boundary matching' density and perform the corresponding micro/macro decomposition:
\be 
\label{rhobar-def}\overline \rho(t,x) =\frac{ \cint{f(t,x,\cdot)}_{V_{-}}}{ \cint{\mathcal E (t,x,.)}_{V_{-}}}, \qquad \cint{\cdot}_{V_{-}}=\int_{V_{-}} \cdot\  d\mu, \qquad  f= \overline \rho \mathcal E + g.
\ee
When $\eps \to 0$, we know that the solution $f$ of \fref{eq-transp} is  (at least formally) close to $\rho \mathcal E$ except in initial or boundary layers. Therefore the 'boundary matching' density
$\overline \rho$ will be close to $\rho$ for small $\eps$.  This  shows that \fref{rhobar-def} still a decomposition of $f$ into an asymptotic part (macro part) and a kinetic part (micro part).
In order to derive the system of equations satisfied by $\overline \rho$ and $g$ from  \fref{eq-transp}, we first integrate \fref{eq-transp} on $V_{-}$ and get
\be
\label{equ-rhobar}
\partial_{t}\overline \rho + \frac{1}{\eps}\cint{v\mathcal E}_{V_{-}}\cdot  \nabla _{x}\overline\rho + \frac{1}{\eps}\cint{v\cdot \nabla _{x}g}_{V_{-}}= \frac{1}{\eps^2}\cint{Lg}_{V_{-}}.
\ee
Substracting this from the equation on $f$, we obtain the equation on $g$:
\be
\label{equ-g1}
\partial_{t}g + \frac{1}{\eps}\left( v\cdot \nabla_{x}g -\cint{v\cdot \nabla_{x}g}_{V_{-}}\right)   + \frac{1}{\eps}\left( v\mathcal E -\cint{v\mathcal E}_{V_{-}}\right)\cdot \nabla_{x}\overline \rho = \frac{1}{\eps^2}\left( Lg - \cint{Lg}_{V_{-}}\right).
\ee
System \fref{equ-rhobar}-\fref{equ-g1} can be replaced by a the more convenient (and still equivalent) system in terms of $\rho =\cint{f} = \int_{V} f d\mu$ and $\ds g= f- \overline \rho \mathcal E$ as follows:
\begin{eqnarray}
\begin{array}{l}
\label{MM-rho-g}
 \partial_{t} \rho +  \frac{1}{\eps}\cint{v\cdot \nabla _{x}g}=0,  \\
\partial_{t}g + \frac{1}{\eps}\left( v\cdot \nabla_{x}g -\cint{v\cdot \nabla_{x}g}_{V_{-}}\right)   + \frac{1}{\eps}\left( v\mathcal E -\cint{v\mathcal E}_{V_{-}}\right)\cdot \nabla_{x}\overline \rho = \frac{1}{\eps^2}\left( Lg - \cint{Lg}_{V_{-}}\right).
\end{array}
\end{eqnarray}
Note that $\overline \rho$ is linked to $\rho$ and $g$ by the relations $ \overline \rho = \rho - \cint{g}= \rho -\cint{g}_{V_{+}}$ and  $f= \rho -\cint{g}_{V_{+}} + g.$

One main interest of this new  micro/macro formulation  \fref{MM-rho-g} for the original kinetic equation  \fref{eq-transp} is the following: the fluxes involved in \fref{MM-rho-g} only concern the quantities $g$ or $\overline \rho$, and not $\rho$. This means that numerical schemes of this formulation would only need the values of $\overline \rho$ and $g$ at the space boundary which are completely known from the original boundary condition  \fref{BC-f} on $f$. 
\section{The case of a fluid scaling}
We now consider non linear  kinetic equations with a fluid scaling
\be
\label{eq-boltz}
   \dt f + v\cdot \grdx f = \frac{1}{\eps} Q(f,f), \qquad t > 0,
  (x,v) \in \Omega \times V, \\   \qquad     f(t=0,x,v) = f_{init}(x,v),
\ee
where $\Omega\subset\RR^d$ and $V=\RR^d$. The collision operator $Q$ is a quadratic operator acting  only on the velocity dependence of the distribution function $f$.  Fundamental examples are the well known Boltzmann kernel
 for rarefied gases and the Fokker-Planck-Landau operators for plasmas. 
In all what follows, we use the notations
\begin{equation} \label{eq-m}
  m(v) = \displaystyle (1, v, |v|^{2}/2)^{T}, 
  \quad \text{ and }\quad \cint g = \int_{\R^{d}} g(v)\ dv
\end{equation}
for any scalar or vector function $g= g(v)$.  It is well known that the Boltzmann  and Landau operators
$Q(f,f)$ have important physical properties  of conservation and entropy: for all $f\geq 0$, we have $\cint {m Q(f,f)}   = 0$ and $\cint {Q(f,f)\log(f)} \leq 0$. It is also well known that the equilibrium functions ($f$ such that $Q(f,f)=0$) are Maxwellians: 
  \begin{equation}\label{Maxwellian2}
    M(U)(v) = \rho(2\pi T)^{-d/2}\exp
    \Big(-|v-u|^{2}/2T \Big),  \qquad U= (\rho, \rho u,  \rho |u|^2/2 + (d/2)\rho T)=\cint{m M(U)} .
  \end{equation}
  To any distribution function $f$, we shall associate its Maxwellian $M=M[f]=M(U)$, that is the fonction of the form
  \fref{Maxwellian2} which has the same moments as $f$:  $U=\cint{mf}= \cint{mM[f]}=\cint{m M(U)}$.
 
When $\eps$ goes to $0$, $f$ approaches a Maxwellian.  At the first order in $\eps$,  the solution $f$ approaches that of the compressible Euler system of gas dynamics. This is formally obtained
by integrating the kinetic equation \fref{eq-boltz} against $m(v)$ and replacing $f$ by its first order approximation: the Maxwellian which has the same first moments as $f$. At the second order in $\eps$, a standard Chapman-Enskog expansion leads to the compressible  Navier-Stokes system.

In order to design numerical schemes which are able to reproduce these asymptotics, a general class of numerical schemes has been introduced in \cite{BLM} and  \cite{Lem-note1} on the basis of micro/macro decompositions. The strategy in \cite{BLM,Lem-note1} consists in decomposing $f$ as $f=M(U) +  g$ and  transform \fref{eq-boltz} into an equivalent system of equations on $M$ and $g$.   Then, discretizations of this formulations lead to a class of AP numerical schemes which are consistent with the fluid limit in the asymptotics $\eps \to 0$.
However, similarly as above for the diffusion limit, the exact space boundary conditions are lost in this decomposition, and one should compute approximate BC for the asymptotic models.  To solve this problem, we introduce here a more suited decomposition matching the boundary.
We first define the new moments $\overline U$ by
$
\overline U(t,x)= \cint{m(v) f(t,x,v)}_{V_{-}}
$
where $V_{-}$ is defined in \fref{VOpm}, and the corresponding equilibrium $\overline M$ is the (uniquely defined, see \cite{Lem-Meh-AP1} for details) Maxwellian such that $\int_{V_{-}} m(v) \overline M dv=\overline U(t,x)$. In some sense, this strategy has some similarities with the half-moment method developed in \cite{dubroca}. Note that $\overline M$ is of the form \fref{Maxwellian2}, i.e. $\overline M=M(\widetilde U)$, but in general $\widetilde U$ is different from $\overline U$.  In fact, since the domain of integration is now $V_{-}$ and not $\RR^d$, the  relation   between  the parameters of  $\overline M$ and $\overline U$ is not explicit in general. We now introduce the following micro/macro decomposition:
\begin{equation}\label{decomp2_f}
  f = \overline M +  g,    \qquad    \int_{V_{-}} m(v) \overline Mdv= \overline U(t,x)= \int_{V_{-}} m(v)f(t,x,v)  dv.
\end{equation}
When $\eps \to 0$,  it is known (at least formally) that the solution $f$ of \fref{eq-boltz} approaches its classical Maxwellian $M$: $ f-M=O(\eps)$. This means that $\int_{V_{-}} m(v)(f-M)= O(\eps)$, and then from \fref{decomp2_f}, we have $\int_{V_{-}} m(v)(\overline M -M)= O(\eps)$. This gives  $d+2$ relations between the two Maxwellians $M$ and $\overline M$, each of these Mawellians being  completely determined by $d+2$ parameters (from \fref{Maxwellian2}). Therefore $\overline M -M = O(\eps)$. In some sense, this justifies the name 'micro/macro' for this new decomposition.

In order to transform \fref{eq-boltz} into a coupled system on $\overline M$ and $g$, we first introduce a definition. For each Maxwellian $\overline M$, we define $\overline \Pi$ as the orthogonal projection on $\mbox{Span}\{\overline M, v\overline M, |v|^2\overline M\}$ in $L^2\left(\overline M^{-1} \un_{V_{-}}\right)$. This projector has an explicit expression, see \cite{Lem-Meh-AP1}.
We now insert the decomposition \fref{decomp2_f} into \fref{eq-boltz}:
$$ \partial_{t} \overline M + v\cdot \nabla_{x} \overline M + \partial_{t} g + v\cdot \nabla_{x}g= \frac{1}{\eps} \left( L_{\overline M} g + Q(g,g)\right), \qquad L_{\overline M} g = 2Q(\overline M,g).$$
Now instead of writing a system on $\overline M$ and $g$, it turns out to be more convenient to rather write a system on the full moments $U=\cint{mf}$ and $g$.  To obtain it, we simply integrate \fref{eq-boltz} against $m(v)$ on the whole velocity domain $V=\RR^d$ for the equation on $U$, and apply $I-\overline \Pi$ to \fref{eq-boltz} to get  the equation on $g$:
\be
\label{MM-fluide}
\begin{array}{l}\partial_{t} U + \nabla_{x}\cdot \left( \int_{\RR^d} v m(v) \overline Mdv\right) + \nabla_{x}\cdot \left( \int_{\RR^d}v m(v)gdv\right) =0, \\
\ds \partial_{t} g + \left(I-\overline \Pi\right)(v\cdot \nabla_{x} \overline M) + \left(I-\overline \Pi\right) \left(v\cdot \nabla_{x} g\right) = \frac{1}{\eps} \left(I-\overline \Pi\right)\left( L_{\overline M} g + Q(g,g)\right).
\end{array}
\ee
Note that, since $f=\overline M +g$, the parameters  $(\overline \rho, \overline u, \overline T)$ of the Maxwellian $\overline M$ (according to definition \fref{Maxwellian2}) are linked to $U$ and $g$ by the simple relation: $ (\overline \rho, \overline \rho \overline u, \overline \rho |\overline u|^2/2 + (d/2) \overline \rho\overline T) = U - \cint{m(v)g}.$

On the other hand, this relation relation cannot determine $\overline M$ at the boundary because $U$ is not known at the boundary. However, by construction, $\overline M$ is completely computed at the boundary from  the relation
$\cint{m(v)\overline M}_{V_{-}}= \cint{m(v)f_{b}}_{V_{-}},$
where $f_{b}$ is the incoming boundary condition \fref{BC-f}.
\section{A numerical test}
\begin{figure}
\begin{center}
\vspace*{-3cm}
\begin{tabular}{@{}c@{}c@{}}
   \hspace*{-.7cm}\includegraphics[width=9cm]{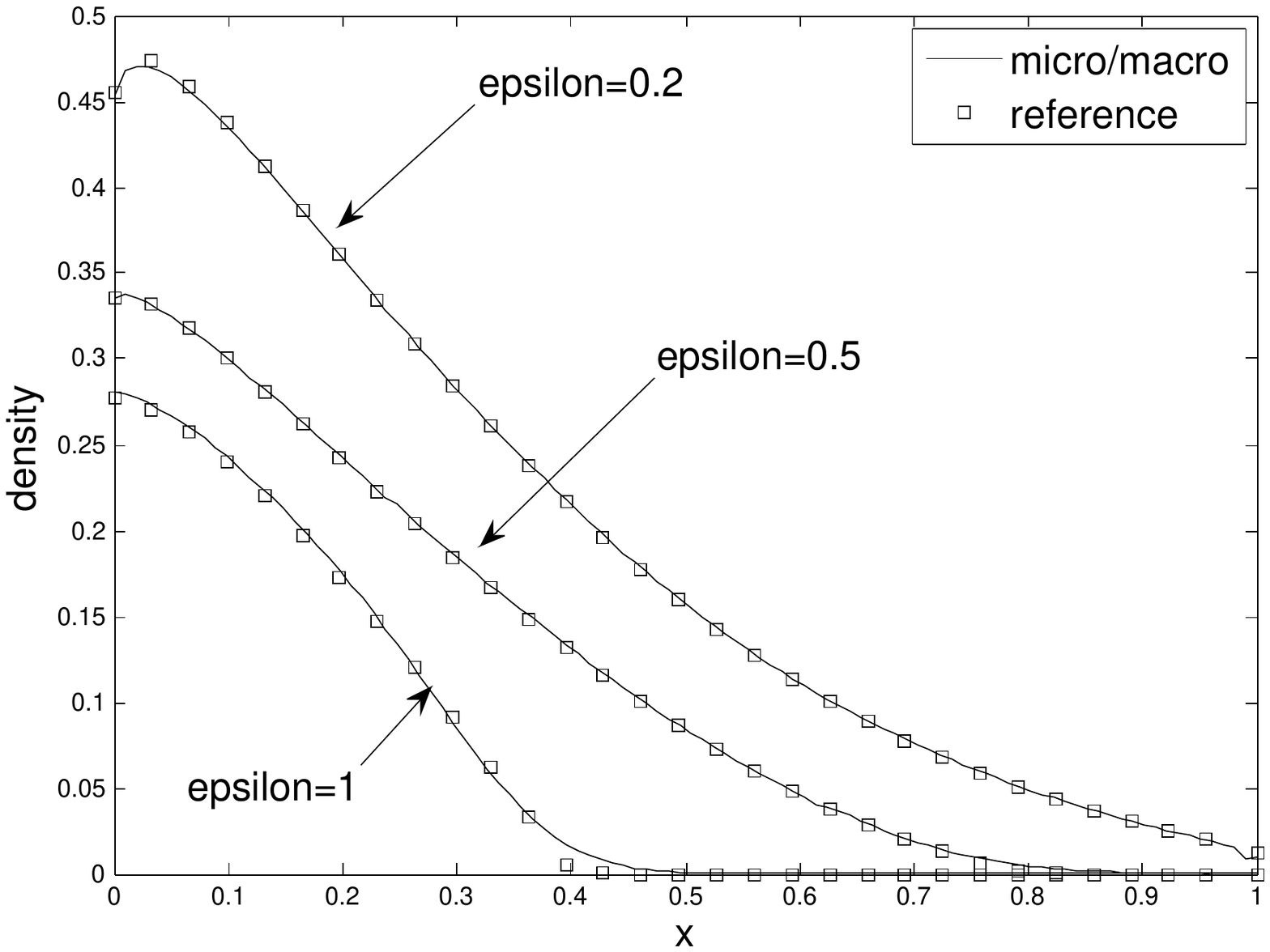}& \hspace*{-.5cm}\includegraphics[width=9cm]{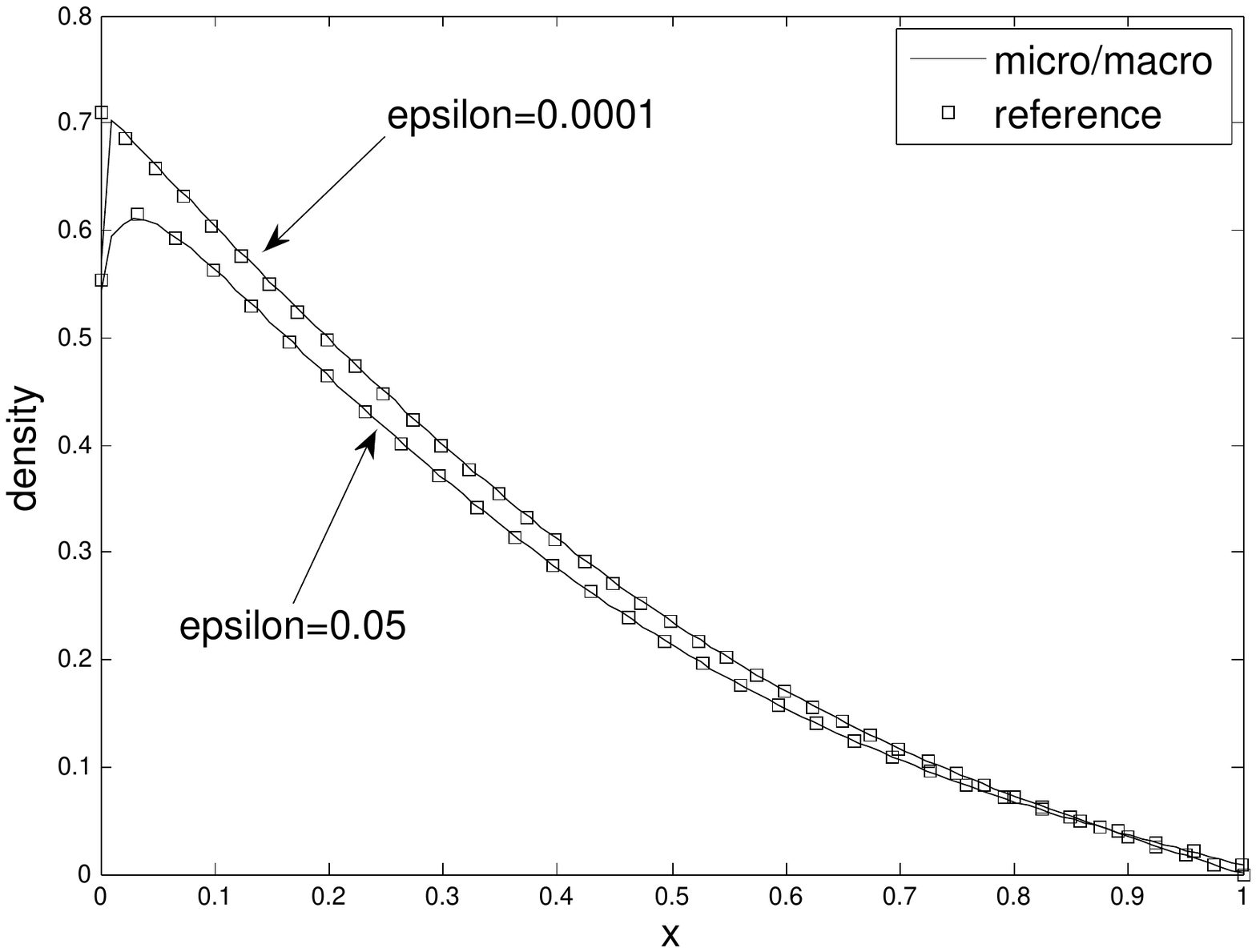} \\
\end{tabular}\vspace*{-3cm}
\caption{Left: numerical solutions in the kinetic and transition regimes, $\eps=1,\,\eps=0.5,\,\eps=0.2$. The density for the micro/macro scheme coincides with the reference calculated with a highly resolved explicit scheme.
Right: numerical solutions in the transition and diffusion regimes. The density for the micro/macro scheme coincides with the reference. For $\eps=0.05$, the reference is calculated with the highly resolved explicit scheme. For $\eps=0.0001$, the reference is calculated by the diffusion equation with the Chandrasekhar value at the boundary. Both figures are plotted at time $t=0.4$ and all micro/macro schemes are used with a uniform grid of $100$ points in $x$.}\label{fig}
\end{center}
\end{figure}

In this Note, we propose to validate the strategy in the case of the diffusion limit only. A more detailed numerical validation will be done in \cite{Lem-Meh-AP1}. We consider the simple framework of one-dimensional spaces $x\in \Omega=[0,1]$ and $v\in V=[-1,1]$ with $d\mu=\frac{1}{2}dv$, and we take $Lf=\frac{1}{2}\int_{-1}^1f dv-f=\rho-f$. The initial data is $f_{init}(x,v)=0$ and the boundary data are $f(t,0,v)=v$ for $v>0$ and $f(t,1,v)=0$ for $v<0$. On Figure \ref{fig}, we plot the densities $\rho$ obtained by the AP micro-macro scheme matching the boundary, for the values $\eps\in\{1,0.5,0.2,0.05,0.0001\}$ and compare them with the reference solutions. For $\eps\in \{1,0.5,0.2,0.05\}$, the reference solutions are computed with a highly resolved explicit scheme for \fref{eq-transp}. For $\eps=0.0001$, the reference solution is computed by the diffusion equation \fref{eq-diff} with the theoretical value of the Dirichlet boundary condition, which is $\rho(0)=0.71043...$ at $x=0$, calculated thanks to the usual Chandrasekhar function, and $\rho(1)=0$ at $x=1$.  The numerical results show that our scheme is AP in the diffusion regime and is able to reproduce efficiently the boundary layer.
The present numerical experiment confirmes that the scheme is 'BAP' in the sense described above: without introducing any artificial boundary condition, we recover both the diffusive regime and the boundary layer.

\bs
\mbox{}

\bs

\vspace*{5mm}

\begin{acknowledgement}
The authors were supported by the french ANR project CBDif.  M. Lemou acknowledges support from the project 'D\'efis \'emergents' funded by the university of Rennes 1. F. M\'ehats acknowledges support from the french ANR project QUATRAIN and from the INRIA project IPSO.
\end{acknowledgement}

\end{document}